\documentclass{amsart} 
\usepackage{amssymb}
\usepackage{amsmath}
\usepackage{amsfonts}

\begin{document}
\newtheorem{theo}{Theorem}[section]
\newtheorem{prop}[theo]{Proposition}
\newtheorem{lemma}[theo]{Lemma}
\newtheorem{exam}[theo]{Example}
\newtheorem{coro}[theo]{Corollary}
\theoremstyle{definition}
\newtheorem{defi}[theo]{Definition}
\newtheorem{rem}[theo]{Remark}


\newcommand{\Bb}{{\bf B}}
\newcommand{\Nb}{{\bf N}}
\newcommand{\Qb}{{\bf Q}}
\newcommand{\Rb}{{\bf R}}
\newcommand{\Zb}{{\bf Z}}
\newcommand{\Ac}{{\mathcal A}}
\newcommand{\Bc}{{\mathcal B}}
\newcommand{\Cc}{{\mathcal C}}
\newcommand{\Dc}{{\mathcal D}}
\newcommand{\Fc}{{\mathcal F}}
\newcommand{\Ic}{{\mathcal I}}
\newcommand{\Jc}{{\mathcal J}}
\newcommand{\Lc}{{\mathcal L}}
\newcommand{\Oc}{{\mathcal O}}
\newcommand{\Pc}{{\mathcal P}}
\newcommand{\Sc}{{\mathcal S}}
\newcommand{\Tc}{{\mathcal T}}
\newcommand{\Uc}{{\mathcal U}}
\newcommand{\Vc}{{\mathcal V}}

\newcommand{\ax}{{\rm ax}}
\newcommand{\Acc}{{\rm Acc}}
\newcommand{\Act}{{\rm Act}}
\newcommand{\ded}{{\rm ded}}
\newcommand{\Gm}{{$\Gamma_0$}}
\newcommand{\ID}{{${\rm ID}_1^i(\Oc)$}}
\newcommand{\PA}{{\rm PA}}
\newcommand{\ACA}{{${\rm ACA}^i$}}
\newcommand{\RefP}{{${\rm Ref}^*({\rm PA}(P))$}}
\newcommand{\RefS}{{${\rm Ref}^*({\rm S}(P))$}}
\newcommand{\Rfn}{{\rm Rfn}}
\newcommand{\tar}{{\rm Tarski}}
\newcommand{\UNFA}{{${\mathcal U}({\rm NFA})$}}

\author{Nik Weaver}

\title [Truth and the liar paradox]
       {Truth and the liar paradox}

\address {Department of Mathematics\\
          Washington University in Saint Louis\\
          Saint Louis, MO 63130}

\email {nweaver@math.wustl.edu}

\date{\em December 28, 2011}

\maketitle


\section{}
The liar paradox is the paradox that arises when we try to assess the
truth value of the sentence ``This sentence is not true.'' We call this
sentence the ``liar sentence''. Assuming it is true leads immediately to
a contradiction, so it must not be true. But that is what it asserts, so
it must be true after all, which is absurd.

The paradox may seem frivolous, but what is worrisome about it is not that
there is any pressing need to assign a truth value to the liar
sentence itself. The broad concern is rather that our inability to resolve
the paradox reveals a fundamental deficit in our understanding of truth
and/or logic. If our grasp of either of these topics is faulty, then the
grand edifice of modern mathematics is called into question. If there is
a specific erroneous step in the paradoxical deduction, and we are unable
to pinpoint it, then this means that some deductive principle we believe
to be valid is not.

That is not to say that there is any likelihood of simple arithmetical
truths ever being overturned, under any circumstances. But mathematics
extends well beyond simple arithmetic, and not all of the principles on
which it is based are so clearly established. Without a complete
understanding of the liar paradox and related problems (Russell's
paradox, etc.), we will never be able to definitively say which of these
principles are genuinely trustworthy. We do have apparently reliable
techniques for avoiding paradoxical reasoning, but these are engineering
solutions. They are not based on a fundamental understanding of the
underlying logical issues.

The liar paradox is strikingly resistant to resolution. One's first reaction
may be that the liar sentence is simply meaningless, or that it does not have
a truth value in a conventional sense. But it is not meaningless to assert
of a meaningless sentence, or one that lacks a truth value, that it is not
true. That kind of assertion is meaningful and correct. Yet this is just
what the liar sentence says of itself, so this naive approach fails
immediately.

The raw intuition that the liar sentence lacks meaningful content has
led to a variety of more sophisticated technical responses, none of which,
however, fares much better than the naive response. The liar sentence is
said to be ``ungrounded'' or ``indeterminate'', both of which would seem
to imply that it is, indeed, correct in saying of itself that it is not
flatly true. Perhaps one wishes to draw a distinction
between grounded and ungrounded truth, but then one has to deal with the
sentence that asserts of itself that it is not true in either sense. At
that point the distinction no longer seems particularly helpful. While it
is certainly clear that a straightforward attempt to determine the truth
value of the liar sentence leads to an infinite regress, it is not so
obvious how that observation helps to resolve the paradox. We can still
talk about such ungrounded sentences, and they, apparently, can talk
about themselves. We might try to say that the truth value of the liar
sentence is unstable, that it should in some way be thought of as flipping
back and forth between truth and falsehood. Very well, then what about the
sentence that asserts of itself that it is continually false, and never true?
It is hard to see how that sentence could ever be true, again leading us back
to a contradiction.

Another class of approaches bring in ideas of ``context'' or ``situation''.
This seems reasonable at first, since most sentences do not have a fixed truth
value, but rather one which varies depending on the context in which they
are asserted or the situation at which they are directed. On the other hand,
it is precisely the assertions of pure logic, into which category the liar
sentence appears to fall, which are not dependent on context, so the relevance
of these kinds of considerations to the liar paradox is questionable. In
any case, we can easily build context-independence into the liar sentence,
by modifying it to explicitly assert of itself that it is not true in any
context. The idea that there could be contexts in which this modified liar
sentence is true is, again, hard to accept.

One may be tempted to try to somehow preempt the sort of objection we have
been making. We could try to say that every sentence has to be interpreted
in some limited context, even the modified liar sentence just mentioned
which purports to refer globally to all contexts. At this point the analysis
becomes farcical: the claim that one cannot refer globally to all contexts
itself refers globally to all contexts and is therefore self-defeating. So
we do not seem to have made much progress.

A completely different type of solution advises that the liar sentence is
both true and not true. But all this does is to eliminate the difficulty by,
in effect, declaring it not to be a difficulty. Surely, to allow a sentence
to be simultaneously true and not true is merely to change the meaning of
the word ``not''. We may well have a trivial resolution of the liar paradox
as it applies to languages in which ``not'' carries some nontraditional
meaning, but obviously this does not address the original problem.

\section{}
It appears that something must be
wrong with our naive conception of truth. This is surprising because it
seems like truth is a definite property of sentences whose
content is exactly captured by Tarski's biconditional \cite{Tar}
$$\mbox{$\ulcorner A\urcorner$ is true}\quad\leftrightarrow\quad
A,\eqno{(*)}$$
where $A$ stands for an arbitrary sentence. We would have to accept that
this intuition is somehow mistaken. But perhaps our intuition about truth
was developed on the basis of unproblematic cases and we are wrong to assume
that it remains valid when applied to exotic (e.g., self-referential)
sentences. Indeed, it is not so clear that the notion of truth is
well-defined in the case of the liar sentence, which strongly suggests
that the naive picture may be misleading, or at least incomplete.

If we decide that our intuition about truth is unreliable, then the paradox
loses much of its force. It merely shows that there can be no predicate
with the properties that we expect a truth predicate to have. Our focus
then shifts to the more straightforwardly technical question, to what
extent is it possible to use ($*$) to define a workable concept of truth?

Since ($*$) tells us, for any sentence $A$, the exact condition for $A$ to
be true, one might think that we could assert it as a definition of truth
for all sentences. But this is impossible, for purely syntactic reasons.
The expression ``For every sentence $A$, $\ulcorner A\urcorner$ is true
if and only if $A$'' is ungrammatical because the variable $A$ functions
ambiguously. In the quantifying phrase ``For every sentence $A$'' it acts
as a reference to an arbitrary sentence, not an assertion of it; but on
the right side of the biconditional it has to be construed as having
assertoric force. The entire expression is ungrammatical in the same way
that ``There exists a sentence $A$ such that $A$'' is ungrammatical. We
should rather say ``There exists a sentence $A$ such that $A$ is true,''
except that this assumes we already have a truth predicate which applies
to arbitrary sentences. In other words, in order to use ($*$) to formulate
a global definition of truth we need some device for asserting a mentioned
sentence, which is just to say that we need to already have a global
truth predicate.

We must understand ($*$) to be not a genuine assertion, but rather a
``schematic'' assertion, something which only becomes a genuine assertion
when some particular sentence is substituted for $A$. Thus, although we cannot
quantify over $A$, we can legitimately use ($*$) to define the truth of
individual sentences by saying things like ``$\,$`Snow is white' is true
if and only if snow is white.'' Modulo finitistic concerns, we could
even conjoin all substitution instances of ($*$), as $A$ ranges over all
the sentences of some language, and thereby produce a truth definition for
that language. But this construction would have to take place within some
(infinitary) metalanguage, and there would always be sentences in that
metalanguage to which the truth predicate just defined did not apply. In
particular, the truth definition would itself be such a sentence, as it
could never appear in one of its own conjuncts.

But even this schematic interpretation does not straightforwardly work for
sentences that themselves refer to truth. Assuming we have abandoned our
naive intuition for truth and are now treating the word ``true'' as initially
undefined, it follows that sentences that explicitly refer to truth are
not initially meaningful and hence cannot be subtituted into ($*$) to
produce a definition of truth. That would be circular, because the
substituted sentence needs to function assertorically on the right
side of ($*$).

Despite this difficulty, there are ways to define truth predicates that
apply in at least some cases to sentences which themselves contain truth
predicates. The most obvious way to do this is, following Tarski, to
construct a well-ordered sequence of truth predicates true${}_\alpha$,
each of which satisfies ($*$) when applied to sentences containing only
truth predicates with lower indices. This resolves the circularity
issue. Another natural approach, due to Kripke \cite{Kri}, is to work
with a single truth predicate but to define its range of application in
a well-ordered sequence of stages. The idea here is that even a partial
definition of truth may be sufficient to give some sentences that contain
the word ``true'' a definite truth value. So every time the range of
application of the truth predicate grows, more sentences are assigned
truth values, and this creates an opportunity to expand the range of
application even further.

\section{}
Both the Tarskian and the Kripkean schemes are natural solutions of
the problem of using ($*$) to define truth. However, neither can be
considered a complete solution. There is a ``revenge'' problem arising
from liar type sentences targeted specifically at the two constructions,
namely the sentences ``This sentence is not true${}_\alpha$ for any
level $\alpha$ of Tarski's construction'' and ``This sentence does not
evaluate as true at any stage of Kripke's construction.'' Provided that
we are working with a language that is sufficiently expressive to allow
us to formulate these assertions, it is easy to see that, indeed, the
first of them is not true${}_\alpha$ for any $\alpha$, and the second
does not evaluate as true at any stage. That is to say, both sentences
can be proven. But the truth predicates we have constructed are evidently
too limited to allow us to say that these sentences are true.

Thus, in both cases there are sentences that we can assert but whose
truth we cannot assert, and this means that a further application of ($*$)
would, in both cases, lead to a truth predicate with broader scope. This
phenomenon is obviously not specific to the two constructions under
consideration. To describe it in general, we might wish to define a
partial truth predicate to be any predicate $T$ such that $T(A) \to A$
holds for all sentences $A$. This expresses the idea that everything
recognized as true really is true, but some true sentences might not
be recognized. We would then say that we have no single notion of
truth; instead, we have a hierarchy of partial truth predicates, and
there is a corresponding hierarchy of liar sentences, each of which
escapes paradox by not falling within the scope of the truth predicate
to which it refers. The lesson of the liar paradox would be that
there is no maximal partial truth predicate.

However, we cannot do this. Our definition of partial truth predicates
is illegitimate because the word ``holds'' functions as a synonym for
``is true''. We might just as well say that $T$ is a partial truth
predicate if $T(A) \to A$ is true for all sentences $A$. In other words,
we need to have a global truth predicate before we can say what constitutes
a partial truth predicate. It really is necessary because both the $T$ on
the left and the $A$ on the right of $T(A) \to A$ must be assertoric, so
without the intervention of a truth predicate they cannot be quantified.
We discussed this kind of syntactic problem earlier.

This syntactic difficulty is genuine. If we allowed ourselves to ignore
it then we would have to also allow a sentence which asserts that no
partial truth predicate holds of itself. This sentence would be paradoxical
in the same way as the original liar sentence since we could immediately
substitute it into ($*$) and obtain a partial truth predicate that applies
to it. It therefore appears that the idea of partial truth is just as
problematic as the idea of truth.

Consequently, although we can see that our ability to go beyond the
Tarskian and Kripkean constructions is a general phenomenon, there is
no obvious way for us to properly express this fact. This state of
affairs is clearly unacceptable. We may have succeeded in convincing
ourselves that our naive intuition for a global notion of truth is
illusory, and this could make the ban on universally quantifying ($*$)
palatable. We can see it as accomplishing the benign task of preventing
us from fomulating a bad definition that affirms a faulty intuition
and has paradoxical consequences. But our corresponding inability to
define the notion of a partial truth predicate is harder to swallow.
It is quite clear that we do have an open-ended ability to use ($*$) to
construct partial truth definitions, and yet somehow we cannot say this.

The same problem was actually present in our initial discussion of ($*$).
When we made the comment that ($*$) could be used schematically to define
the truth of any sentence, we did not notice that this statement, too,
cannot be grammatically formulated unless a global truth predicate is
already available. We cannot say ``For any meaningful sentence $A$, ($*$)
can be used to define a predicate $T$ such that $T(A) \leftrightarrow A$.''
Proper syntax demands something like ``For any meaningful sentence $A$, ($*$)
can be used to define a predicate $T$ such that $T(A) \leftrightarrow A$ is
true'' or ``such that $T(A)$ is true if and only if $A$ is true'', so again,
we need to already have a global truth predicate.

\section{}
What has been missing from this discussion is an explicit acknowledgement
of the constructive quality of Tarski's biconditional. We alluded to it
just now in our comment that we have an open-ended ability to construct
partial truth definitions. It is inherent in the nature of ($*$) both
that we have a global ability to use it to define truth in limited settings,
and that we can always go beyond any such setting to produce a more inclusive
definition. We are having difficulty precisely expressing this idea because
we are trying to express it in classical terms.

Fortunately, there is a well-developed theoretical apparatus for dealing
with constructive phenomena. We have a form of reasoning, intuitionistic
logic, which is suited to this kind of setting, and we have an interpretation
of the logical constants, the BHK or ``proof'' interpretation, which
supports this form of reasoning. The identification of provability as the
key primitive concept in terms of which constructive reasoning is based
is crucial because it provides us with the linguistic tool that
we need to handle use/mention problems of the type that arose in our attempts
to define truth and partial truth. In particular, the constructivist
solution to the problem of defining truth is to simply equate it
with provability. Thus, letting $\Box A$ stand for ``$A$ is provable'',
i.e., there exists a proof of $A$, we can say
$$\mbox{$A$ is constructively true}
\quad\leftrightarrow\quad \Box A.\eqno{(**)}$$
Observe that neither appearance of $A$ in ($**$) is assertoric, so there
is no syntactic obstruction to quantifying over $A$ in this expression.
We are free to say that for any sentence $A$, $A$ is constructively true
if and only if $A$ is provable. Thus the constructive notion of truth is
genuinely global in scope.

The other problematic aspect of ($*$), that it becomes circular when
applied to sentences containing the word ``true'', also does not translate
to ($**$). We take provability to be an objective
notion that exists independently of our characterization of it, so that
there is nothing ill-defined about the provability of
sentences which themselves refer to provability.

The difficulty we encountered in trying to define partial truth predicates
can be avoided in a similar way. We can define a predicate $T$ to be a
partial truth predicate if for any sentence $A$ the implication $T(A) \to A$
is provable. This finesses the syntactic problem discussed earlier and
provides us with a global notion of partial truth. We might also say that a
sentence $A$ falls within the scope of a partial truth predicate $T$ if
the biconditional $T(A) \leftrightarrow A$ is provable, and that
one partial truth predicate, $T_1$, is subordinate to another, $T_2$, if
for every sentence $A$ the implication $T_1(A) \to T_2(A)$ is provable.
This enables us to make the comment that we tried to make earlier, about
there being a hierarchy of partial truth predicates and a corresponding
hierarchy of liar sentences.

Note that under the above definition of partial truth there is no requirement
that the predicate $T$ must itself be constructive. To the contrary, we can
define partial truth predicates in wide generality using either Tarskian or
Kripkean techniques. To give one example, we can define Goldbach's conjecture
to be true if and only if every even number greater than two is a sum of two
primes. Whether Goldbach's conjecture has a constructive truth value is
immaterial here; what matters is that under this definition Goldbach's
conjecture and the truth of Goldbach's conjecture are provably equivalent.
Thus, although the notion of what constitutes a partial truth
predicate is constructive, partial truth predicates may themselves
be highly nonconstructive.

A variety of notions of truth are available. There is a global constructive
truth predicate and there is a hierarchy of classical truth predicates with
limited scope. What we cannot have is a global classical truth predicate.

\section{}
Bringing constructive ideas into play has given us the ability to
globally discuss partial truth predicates, but by introducing a new,
constructive notion of truth it also exposes us to a new version of
the liar paradox. What are we to make of the ``provable liar'' sentence
that asserts of itself that it is not provable, or the sentence that
asserts that no partial truth predicate provably holds of itself?

Since the notion of constructive truth is global in scope, the provable
liar sentence cannot be handled in the same way as the classical liar
sentence. On the other hand, a paradox is less immediate here because
the validity of Tarski's biconditional is not directly evident in the
constructive setting. We cannot just assume that asserting $A$ is
equivalent to asserting that $A$ is provable. There is clearly some
relationship between the two statements, but we have to determine
what that relationship is.

At first sight it appears that asserting $A$ is provable is stronger
than merely asserting $A$. In one direction, we have no reason to expect
that every classically true statement must be provable. But in the other
direction, it seems clear that valid reasoning must always lead to a true
conclusion. If there is any doubt about this, surely we can simply build
it into our notion of what counts as a valid proof. Thus, we can apparently
affirm $\Box A \to A$ but not its converse.

However, this reasoning assumes that we are using classical logic and
that a classical notion of truth is available. If we are working in a
setting in which atomic sentences might not have definite truth values,
or we lack a classical truth predicate, then the conclusion we just reached
must be reconsidered. First, there could be a problem with the argument
that we are free to fiat that valid proofs respect classical truth; second,
if a truth table definition of implication is unavailable, then we have
to switch to a constructive interpretation of implication, which
changes the nature of the question.

So before we go any further with the provable liar sentence we have to
decide whether the concepts present in it should be understood classically
or constructively. This comes down to the question of whether we have the
right to assume that $\Box A$ has a definite truth value for every
sentence $A$.

In the constructivist literature one sometimes sees the claim that every
assertion
of the form ``$p$ proves $A$'' not only has a definite truth value but is
even, in principle, decidably true or false. However, this may be another
case where our intuition is developed on the basis of unproblematic cases
and becomes unreliable when we have to deal with self-referential phenomena.
There could be situations where in the course of assessing the validity of
some alleged proof we need to first assess the validity of some other
alleged proof, and this creates the possibility for a vicious circle.

For instance, let $A$ be the sentence ``The sentence $p$ is not a proof of
this sentence'' and let $p$ be the sentence ``By inspection, this sentence
is not a proof of the sentence $A$.'' If $p$ is a proof of $A$ then a
falsehood is provable, so we surely want to say that $p$ does not prove $A$.
However, decidability of the proof relation implies that if $p$ does not
prove $A$, then this can be seen by inspection. If we accept this principle
and we agree that $p$ does not prove $A$ then we should also agree that the
argument ``By inspection, $p$ is not a proof of $A$'' is a valid proof of the
fact that $p$ is not a proof of $A$, which is just to say that we have
to agree that $p$ is a proof of $A$. So there is a paradox.

A contradiction can be avoided by adopting the position that $p$ is not a
proof of $A$, but that this cannot be seen by inspection.
But this does not really help, because if some argument has to be
made to show that $p$ is not a proof of $A$ then we ought to be able to
incorporate that argument into $p$. For instance, we could replace $p$
with the sentence ``If this sentence were a proof of $A$ then a falsehood
would be provable, so this sentence is not a proof of $A$'' or something
of that sort. A moment's thought shows that the even more extreme position
that $p$ is not a proof of $A$, but there is no way to see this, is
straightforwardly self-defeating. (How could we know this?)

The explicit labelling of $A$ and $p$ is not necessary either. We could
take $A$ to be, say, the sentence ``There is no proof of this sentence that
is fewer than 1,000 characters long'' and $p$ to be the sentence ``By
inspection, no string of fewer than 1,000 characters is a proof of the
sentence `There is no proof of this sentence that is fewer than 1,000
characters long'.'' (Or: if any such string were a proof of that sentence
then a falsehood would be provable, so no such string is a proof of
that sentence.)

What the preceding should show is that it is not so clear that every
statement of the form ``$p$ proves $A$'' has a decidable (or even definite)
truth value. We must be careful to distinguish decidability of the proof
relation from the weaker claim that any proof can be recognized to be a
proof. The latter clearly is inherent in our notion of what proofs are;
if a proof is an argument which is completely persuasive to a rational
mind, then certainly, in order for something to count as a proof its
validity must in principle be recognizable. But to infer decidability of
the proof relation we would have to be able to also affirm that anything
which is not a proof can be recognized not to be a proof. (Any Turing
machine which halts on null input can be recognized to halt on null input, but
that does not mean we can always decide whether a given Turing machine halts on
null input.) This direction is not obviously inherent in our notion of proof.

We conclude that reasoning classically about the general notion of
provability is unjustified. We have to reason constructively.

\section{}
Since the provable liar sentence explicitly involves the general
notion of provability, it must be understood constructively. We now need
to determine the relationship between $A$ and $\Box A$ under constructive
logic. The relevant question is whether it is always possible
to convert a proof of $A$ into a proof of $\Box A$, and vice versa.

First, we need to clarify what would constitute a proof of $\Box A$.
By definition, to say that $A$ is provable is to say that there exists
a proof of $A$. So, letting $p\,\vdash A$ stand for ``$p$ proves $A$'',
a proof of $\Box A$ is a proof of $(\exists p)(p\,\vdash A)$. But if we
are reasoning constructively, any existence proof should in principle
provide an explicit instance. So $q$ is a proof of $\Box A$ if and only
if, for some $p$, $q$ proves $p\, \vdash A$.

Consider the implication $A \to \Box A$. To prove this implication we
must show how to convert any proof of $A$ into a proof of $\Box A$.
But when we are compelled by a proof $p$ to affirm an assertion $A$,
we are simultaneously compelled to affirm the assertion that $p$
proves $A$. (If we do not realize that we have proven $A$, then we
have not proven $A$.) Thus, whenever $p$ proves $A$ it also proves
$p\, \vdash A$. So any proof of $A$ is also a proof of $\Box A$, and
this shows that the law $A \to \Box A$ is constructively valid.

Now consider the converse implication $\Box A \to A$. Here we must
show how to convert any proof of $\Box A$ into a proof of $A$. The
obvious procedure in this case is, given $p$ and $q$ such that $q$
proves that $p$ proves $A$, to discard $q$ and return $p$. This
procedure succeeds if, whenever $q$ proves that $p$ proves $A$, $p$
actually does prove $A$. Informally, we need to know that proofs
are reliable. But that is just what we are trying to verify: in
order to establish that $\Box A$ implies $A$, we already need to
know that $\Box(p\, \vdash A)$ implies $p\, \vdash A$. In contrast
to the $A \to \Box A$ direction, in this direction the special case
of statements of the form $p\, \vdash A$ is no more evident than the
general case. The obvious attempt to justify the law $\Box A \to A$
is therefore circular.

In any instance where we have actually proven $\Box A$, i.e., we have
proven that some $p$ proves $A$, we ought to be willing to accept $p$
as a proof of $A$ and therefore infer $A$. Making this deduction requires
only that we accept the reliability of the proof we have just given, not
the global reliability of all proofs. So the passage from $\Box A$ to $A$
manifests as a deduction rule rather than an implication.

Can we build reliability into our notion of valid proof? That is, can we
modify our concept of proof so as to explicitly include the requirement
that whenever $p$ proves $A$, $A$ is true? In order to modify our notion
of proof in this way we would need to use a classical truth predicate. But
until the meaning of provability is settled, a classical truth predicate
is not available for sentences like the provable liar which explicitly
refer to the notion of provability. We cannot use ($*$) to define truth
for any sentence whose meaning is not settled, and we have not settled
the meanings of sentences which involve the notion of provability until
we have specified our notion of  proof. So a classical truth predicate for
such sentences is not available for us to use when we are at the stage of
specifying what counts as a proof. Therefore, in
this setting we cannot incorporate the requirement that
only true sentences are provable into our concept of proof.

All our attempts to justify the law $\Box A \to A$ under a constructive
interpretation of implication are ultimately circular. So we cannot
assume this law is generally valid. Only the direction $A \to \Box A$
can be universally asserted.

It may be helpful here to note that a similar restriction appears when we
follow Kripke's prescription for using ($*$) to generate a self-applicative
truth predicate. At any stage of the construction we are free to affirm that
the currently available predicate is sound, i.e., satisfies the law
$T(A) \to A$. That is, we can prove this implication for any
sentence $A$ which lies within the current scope of $T$. However, when
extending the current predicate we do not
assume that subsequent stages will be sound. That would be patently
circular, and it would generate a contradiction in the following
way. If we assume that all future extensions of the truth predicate we
are constructing will be sound, then a short deduction shows that the
liar sentence which states that it does not evaluate as true at any stage
will, indeed, not evaluate as true at any stage. Thus, the soundness of all
future stages entails this liar sentence, and this is all we need to know
in order to be justified
in designating this liar sentence as true at the next stage. In other words,
it is only the possibility that the liar sentence might be recognized as true
at some future stage that prevents us from recognizing it as true now. The
lesson is that it is appropriate to affirm the soundness of earlier stages,
but wrong and dangerous to assume the soundness of later stages. We may
anticipate that future stages will be sound, but we cannot employ this
as a premise when generating the next stage.

The analogy is that disallowing $\Box A \to A$ as a proof principle
involves a similar restriction. Just as we cannot employ the premise
that a not yet fully defined truth predicate is globally sound in order
to establish that a given sentence is true, we cannot employ the premise
that proofs are globally sound in order to establish that a given proof
is valid. Thus, since the justification of $\Box A \to A$ hinges
on the global soundness of all proofs, this justification fails.

\section{}
We are now in a position to determine what can be said about the provable
liar sentence. This sentence asserts that it is not provable, which we can
express symbolically as $L \equiv \neg \Box L$. Now assuming $L$, we can
immediately deduce $\neg \Box L$. But we can also deduce $\Box L$ by
using the general law $A \to \Box A$. Combining these statements, we
get that $L$ entails a contradiction, which means that we have proven
$\neg L$. It follows that we can also prove $\neg\neg\Box L$.

But this does not directly lead to any paradox. If we could infer
$\neg \Box L$ from $\neg L$, then the two conclusions we just reached
would yield a contradiction. But the implication goes in the other
direction: in general we have $\neg \Box A \to \neg A$, not conversely.

A minor variation on the provable liar sentence is the sentence which
asserts that its negation is provable.
Here the formalization is $L' \equiv \Box(\neg L')$,
and we reach slightly different conclusions. The results in this case
are $\neg\neg L'$ and $L' \to \Box \bot$, where $\bot$
represents falsehood. Since $\neg A$ equals $A \to \bot$ by
definition, we may say that the weaker statement $A \to \Box \bot$
represents ``weak'' falsity of $A$. Thus the conclusions we reach
are that $\neg L'$ is false and $L'$ is weakly false.

In order to satisfy ourselves about the consistency of this kind of
reasoning, we can set up a formal constructive self-referential
propositional calculus. This is done in \cite{Wea}, and we show
there that the system is consistent. Thus we are able to reach the
substantive conclusions about $L$ and $L'$ mentioned above without
producing a contradiction.

The analysis of the sentence which asserts that no partial truth
predicate provably holds of itself is similar. Again, the absence
of the law $\Box A \to A$ blocks the paradox. In the case of the
provability paradox discussed in Section 4, the result is that
``$p$ is a proof of $A$'' is weakly false and ``$p$ is not a proof of
$A$'' is false.

Earlier we asserted that there can be no global classical truth predicate.
The way to make this idea precise is to assume that there is a predicate
$T$ such that the biconditional $T(A) \leftrightarrow A$ is provable
for all sentences $A$; then taking $A$ to be the corresponding liar
sentence allows us to deduce $\Box \bot$. So we conclude that the
existence of such a predicate is weakly false.

\section{}
To summarize: the classical liar paradox is vacuous because there
is no global classical truth predicate. We are only able to define
predicates which satisfy Tarski's biconditional in limited settings.
But the global assertion that a truth definition can be given for any
meaningful sentence cannot be stated classically, only constructively.

In the constructive setting truth is equated with provability, and this
is a global notion. Therefore we do have a genuine constructive liar
sentence. But here there is no paradox because neither the law of excluded
middle nor the law $\Box A \to A$ can be universally affirmed.

Attempts to deal with the liar paradox typically propose some sort of
restriction on meaningful or acceptable sentences and then observe that
the liar sentence violates it. But the justification of the restriction,
or even its mere expression, invariably violates the very same principle,
making the whole account self-defeating. We cannot reject all
ungrounded assertions because the sentence which expresses this prohibition
is itself ungrounded. We cannot stratify sentences into types because the
assertion that every sentence has a type does not itself have a type. We
cannot demand that every sentence be restricted to some local context
because that demand itself would have to be restricted to some local
context. And so on.

This is just a consequence of treating a constructive phenomenon as if it
were classical. Indeed, no satisfactory classical conclusion can be drawn
about the provable liar sentence: coming to a firm conclusion about its
truth or falsehood would spell disaster either way. The way forward
is to avoid contradiction by leaving open the possiblity that the liar
sentence might be recognized as having a definite truth value. This
outcome would be disastrous, but that does not license us to
deny its possibility. And that is good, because not denying this
possibility is exactly what blocks the paradox.



\bigskip
\begin{thebibliography}{aaaaaaaa}

\bibitem{Kri}
S.\ Kripke, Outline of a Theory of Truth, {\it Journal of Philosophy \bf 72}
(1975), 690-716.

\bibitem{Tar}
A.\ Tarski, The semantic conception of truth, {\it Philosophy and
Phenomenlological Research \bf 4} (1944), 13-47.

\bibitem {Wea}
N.\ Weaver, The semantic conception of proof, manuscript.\footnote{See
http://www.math.wustl.edu/$\sim$nweaver/conceptualism.html}

\end{thebibliography}
\end{document}